\documentclass[11pt]{amsart}

\usepackage{amsmath,amssymb,amsthm}
\usepackage{mathrsfs}
\usepackage{tikz}
\usepackage[margin=1.2in]{geometry}
\usepackage{url}
\newenvironment{numthm}[2]%
  {\par\medskip\noindent{\scshape #1~#2.}\ \itshape\ignorespaces}%
  {\par\medskip}
\newenvironment{numres}[2]%
  {\par\medskip\noindent{\scshape #1~#2.}\ \ignorespaces}%
  {\par\medskip}
\newenvironment{rem}[1][Remark]%
  {\par\smallskip\noindent\textit{#1.}\ \ignorespaces}%
  {\par\smallskip}

\newcommand{\On}{\mathcal{O}_n}
\newcommand{\Cn}{\mathcal{C}_n}
\newcommand{\E}{\mathbb{E}}
\newcommand{\Var}{\operatorname{Var}}
\newcommand{\rank}{\operatorname{rank}}
\newcommand{\im}{\operatorname{im}}

\title{Combinatorial results for semigroups of order-preserving mappings}
\author{Peter M. Higgins}
\address{Department of Mathematics, University of Essex, Colchester, UK}

\thanks{This is the author's version of the article published in
\emph{Math.\ Proc.\ Camb.\ Phil.\ Soc.} \textbf{113} (1993), 281--296,
\url{https://doi.org/10.1017/S0305004100075964}.
Received 12 March 1992; revised 31 July 1992.}


\begin{document}

\begin{abstract}
We investigate combinatorial properties of the semigroup $\mathcal{O}_n$ of all
order-preserving mappings on $X_n=\{1,2,\dots,n\}$ and of its subsemigroup
$\mathcal{C}_n$ of all decreasing order-preserving mappings. We show that
$|\mathcal{C}_n|$ is the $n$th Catalan number and obtain a new proof of the
standard Catalan identity. Using ballot arguments and the reflection principle
we find the probability that a random member of $\mathcal{O}_n$ is
indecomposable, and via the Wallis product we show that the mean number of
orbits of a random $\alpha\in\mathcal{O}_n$ approaches $\frac12\sqrt{n\pi}$. By
contrast, the mean number of orbits of a random $\alpha\in\mathcal{C}_n$ never
exceeds $3$; we find the distribution of the orbit number explicitly in both
cases and obtain the limiting distribution for $\mathcal{C}_n$.
\end{abstract}

\maketitle

\section{Introduction}

Consider the finite set $X_n=\{1,2,\dots,n\}$ ordered in the standard way. Let
$T_n$ denote the full transformation semigroup on $X_n$, that is, the semigroup
of all mappings $\alpha\colon X_n\to X_n$ under composition. We shall call
$\alpha$ \emph{order-preserving} if $i\le j$ implies $i\alpha\le j\alpha$ for
$i,j\in X_n$, and $\alpha$ is \emph{decreasing} if $i\alpha\le i$ for all
$i\in X_n$. This paper investigates combinatorial properties of the semigroup
$\On$ of all order-preserving mappings on $X_n$, and of its subsemigroup $\Cn$,
which consists of all decreasing and order-preserving mappings.

Combinatorial properties of $T_n$ have been studied over a long period under
various guises and many interesting and delightful results have emerged (see,
e.g.\ \cite{H90,H92}). Two papers concerned with combinatorial results of $\On$
are \cite{Howie} and \cite{GH}, while little seems to have been written on
$\Cn$. The inverse semigroup of partial one-to-one order-preserving maps has
been studied by Garba \cite{Garba}. Nevertheless both $\On$ and $\Cn$ have
arisen naturally in language theory: any $\mathscr{J}$-trivial finite semigroup
divides some $\Cn$ (see \cite{Pin}) while a description of the class of
semigroups that are divisors of some $\On$ remains an interesting open problem.
It comes as a pleasant surprise that many classical results from elementary
analysis and combinatorics emerge as keys to our combinatorial problems. For
instance it transpires that $|\Cn|=C_n$, the $n$th Catalan number, and indeed
the standard Catalan identity (Result~1.4) follows easily from a natural
partition of $\Cn$. Many of our counting tasks can be formulated as certain
`ballot' problems and thus may be solved using Andr\'e's Reflection Principle
(Result~1.3). From single-variable calculus the Wallis Product arises during the
calculation of the mean number of orbits of a randomly chosen $\alpha\in\On$ and
allows us to prove that this number approaches $\frac12\sqrt{n\pi}$ for large
$n$. By way of contrast, the mean number of orbits of $\alpha\in\Cn$ never
exceeds $3$: we find the distribution of the orbit number explicitly in this
case and obtain the limiting distribution.

We list some standard combinatorial results germane to our purposes. We shall
write the binomial coefficient $\binom{n}{r}$ also as $C(n,r)$.

\begin{numres}{Result}{1.1}
\begin{enumerate}
\item[(i)] $\sum_{k=0}^{n} C(n,k)\,C(m,k)=C(n+m,m)$ for $n\le m$.
\item[(ii)] $\sum_{k=0}^{n} C(2k,k)\,C(2n-2k,n-k)=4^n$.
\end{enumerate}
\end{numres}

To see (i), consider a collection of $n$ red and $m$ blue labelled balls. The
number of ways of choosing $m$ balls from this collection can be thought of as
the sum, from $k=0$ to $n$, of the number of ways of choosing $k$ blue balls and
$m-k$ red balls, from which the identity follows.

To prove (ii) one can use the factorization
$(1-4x)^{-1}=(1-4x)^{-\frac12}\cdot(1-4x)^{-\frac12}$. The coefficient of $x^n$
in the binomial expansion of $(1-4x)^{-\frac12}$ is $C(2n,n)$, while the
corresponding coefficient in the geometric series for $(1-4x)^{-1}$ is $4^n$.
The coefficient in the product is then the required summation.

As explained in \cite{Howie}, the members of $\On$ are in a natural one-to-one
correspondence with the sequences of $n$ $A$'s and $n$ $B$'s ending in $B$:
given such a sequence the corresponding mapping $\alpha\in\On$ is defined by the
rule $i\alpha=j$, where $j$ is the number of the next $B$ in the sequence
following the $i$th $A$. From this we see

\begin{numres}{Result}{1.2}
$|\On|=C(2n-1,n-1)$.
\end{numres}

Deleting the final $B$ in such a sequence gives us an arbitrary sequence of $n$
$A$'s and $(n-1)$ $B$'s which can be thought of as the result of counting a
ballot in which the two candidates, $A$ and $B$, receive $n$ and $n-1$ votes
respectively. Such a count can be represented as a lattice path from $(0,0)$ to
$(2n-1,1)$, where each point of the count corresponds to a lattice point
$(p,q)$, and the following point on the path is $(p+1,q\pm1)$ according as $A$
or $B$ polls the next vote. The lattice point $(p,q)$ then signifies a lead of
$q$ votes for $A$ after $p$ votes have been counted.

\begin{numthm}{Result}{1.3 \textup{(the reflection principle)}}
Let $a$, $b$ and $c$ be positive integers. Then the number of lattice paths from
$(0,a)$ to $(b,c)$ which touch or cross the $x$-axis equals the total number of
lattice paths from $(0,-a)$ to $(b,c)$.
\end{numthm}

\begin{figure}[ht]
\centering
\begin{tikzpicture}[scale=1.0,
    solid dot/.style={circle,fill=black,inner sep=0pt,minimum size=5pt},
    open dot/.style={circle,draw=black,fill=white,inner sep=0pt,minimum size=5pt}]
  \draw[thick] (0,-2.6) -- (0,2.6);
  \draw[thick,->] (0,0) -- (10.6,0) node[below] {$x$};
  \draw[thick] (0,1) -- (1,2) -- (2,1) -- (3,0) -- (4,-0.8);
  \foreach \p in {(0,1),(1,2),(2,1),(3,0),(4,-0.8)} \node[solid dot] at \p {};
  \draw[dashed] (0,-1) -- (1,-2) -- (2,-1) -- (3,0);
  \node[open dot] at (1,-2) {};
  \node[open dot] at (2,-1) {};
  \node[open dot] at (0,-1) {};
  \node at (5.6,-0.35) {$\bullet\ \ \bullet\ \ \bullet$};
  \draw[thick] (7.4,1) -- (8.6,0) -- (9.8,1);
  \foreach \p in {(7.4,1),(8.6,0),(9.8,1)} \node[solid dot] at \p {};
  \node[left] at (-0.1,1) {$(0,a)$};
  \node[left] at (-0.1,-1) {$(0,-a)$};
  \node[right] at (9.9,1) {$(b,c)$};
\end{tikzpicture}
\caption{The reflection principle.}
\label{fig:reflection}
\end{figure}
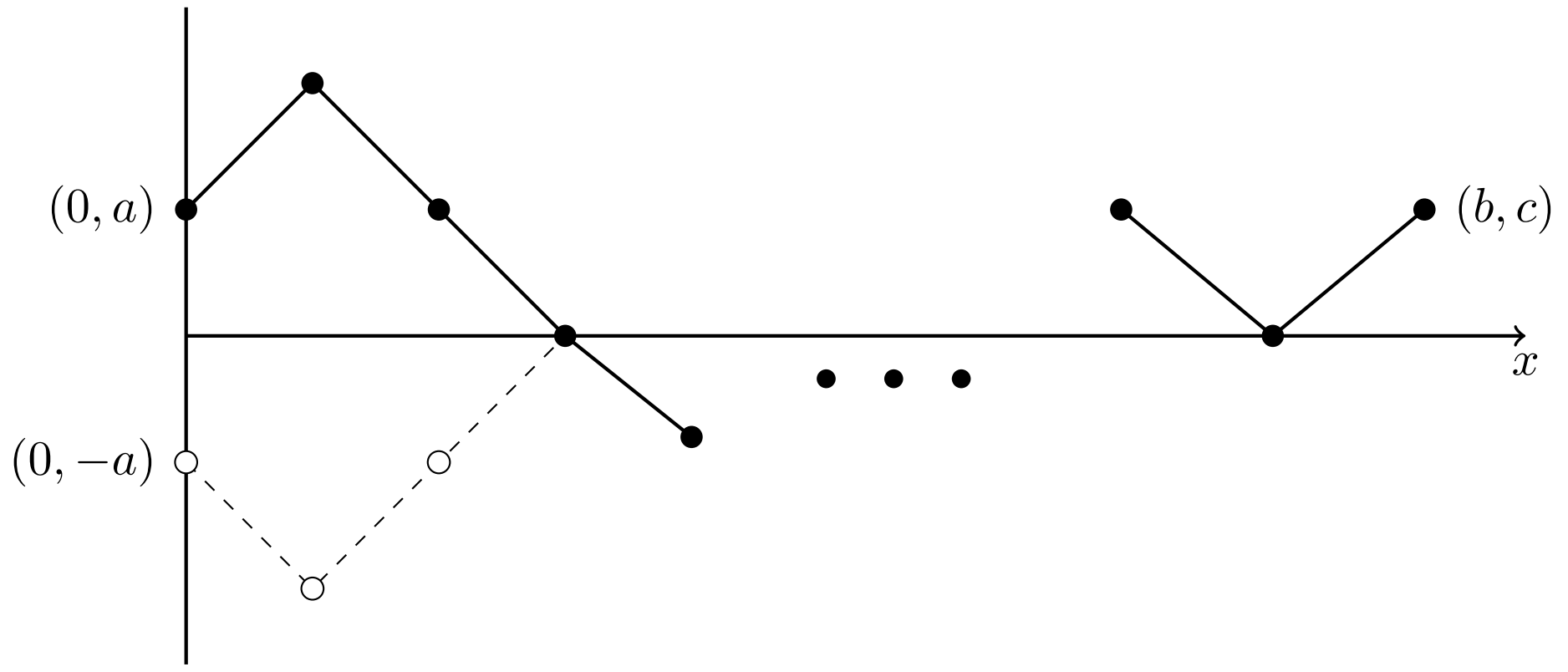

The principle follows from the observation that there is a bijection between the
two types of path defined by reflecting in the $x$-axis the initial segment of
any path of the first type up to the point where it first meets the axis (see
Figure~\ref{fig:reflection}).

The bold path represents a typical lattice path from $(0,a)$ to $(b,c)$ which
first meets the axis as indicated. The dotted path represents the reflection of
the initial segment which, together with the remainder of the path, yields a
lattice path from $(0,-a)$ to $(b,c)$. In general there will be paths of either
type if and only if there is some positive integer $d$ (representing the number
of down-directed edges of the path from $(0,-a)$) such that $2d+a+c=b$, that is,
if and only if $b-(a+c)$ is an even positive integer. The number of such lattice
paths is then $C(b,d)$, corresponding to the number of choices for selecting the
positions of the $d$ down-edges amongst the $b$ edges of the lattice path.

The $n$th Catalan number $C_n$ is $C(2n,n)/(n+1)=C(2n,n-1)/n$. Historically,
$C_n$ was defined as the number of distinct bracketings of the $(n+1)$-fold
product $x_1x_2\cdots x_{n+1}$. Equally, $C_n$ may be defined as the number of
binary trees with $n$ source nodes, and as the number of ways of triangulating a
convex $(n+2)$-gon by means of non-intersecting diagonals (see \cite{HP}). The
amenability of Catalan numbers to manipulation stems from the following
convolution-style identity.

\begin{numres}{Result}{1.4}
$\sum_{k=1}^{n} C_{k-1}C_{n-k}=C_n$.
\end{numres}

This result can be verified by partition arguments appropriate to the type of
collection used to define $C_n$: we provide a new argument based on partitioning
$\Cn$ in Section~3.

Section~2 deals with combinatorial aspects of $\On$ while corresponding results
for $\Cn$ appear in Section~3. Some results for $T_n$ are stated along the way
for the sake of comparison. The content of Section~2 was in part the subject of
the author's talk at the Lisbon Conference on Lattices, Semigroups and Universal
Algebra, and appears in the proceedings of that meeting \cite{H90}.

\begin{rem}[Notation]
We shall write $\mu_X$, $\sigma_X$, $\sigma_X^2$ for the mean, standard
deviation, and variance respectively of a random variable $X$ (we shall also
write $\Var(X)$ for the latter quantity). Each $\alpha\in T_n$ can be pictured
as a digraph on $n$ vertices with $ij$ an arc of $\alpha$ if $i\alpha=j$. Each
component of such a digraph is \emph{functional}, meaning that it consists of a
unique cycle, together with a number of trees rooted around the points of the
cycle. The components of the digraph of $\alpha$ correspond to the orbits of the
mapping $\alpha$, where we say that the \emph{orbit} of $i\in X_n$ is
$\{j\in X_n: i\alpha^r=j\alpha^s \text{ for some positive integers } r,s\}$.
The arrows within a tree rooted on the cycle of some component of $\alpha$ point
along the path towards the root of the tree. Hence if we follow the convention
that the cycles of $\alpha$ are directed anti-clockwise, then the arrows may be
deleted from the picture of $\alpha$ without loss of information, provided that
one-cycles (which correspond to fixed points) are shaded to avoid ambiguity. For
further background see \cite{H92}.
\end{rem}

\begin{numthm}{Proposition}{1.5}
The cycles of the components of $\alpha\in\On$ each consist of a unique fixed
point. Each orbit of $\alpha\in\On$ is convex in the ordered set $X_n$.
\end{numthm}

\begin{proof}
Suppose that $i,j\in X_n$ are distinct members of some cycle of a component of
$\alpha$, so that $j=i\alpha^r$, $i=j\alpha^s$ for some positive integers $r$
and $s$. Suppose that $i<i\alpha$. We then obtain an increasing sequence
$i,i\alpha,i\alpha^2,\dots,i\alpha^r=j,i\alpha^{r+1},\dots,i\alpha^{r+s}$. But
$i\alpha^{r+s}=j\alpha^s=i$, giving equality throughout, and yielding the
contradiction that $i=j$. The same contradiction results if $i>i\alpha$, so that
all cycles of the digraph of $\alpha\in\On$ are trivial.

Next let $C$ be an orbit of $\alpha\in\On$ with (unique) fixed point $p$ say,
and let $c$ and $d$ be the respective minimum and maximum members of $C$. Let
$r,s$ be the respective least non-negative integers such that
$c\alpha^r=d\alpha^s=p$. Then
\[
c<c\alpha<c\alpha^2<\cdots<c\alpha^r=p=d\alpha^s<d\alpha^{s-1}<\cdots<d\alpha<d.
\]
Let $k\in X_n$ lie between $c$ and $d$. Then either
$c\alpha^m\le k\le c\alpha^{m+1}$ for some $0\le m\le r-1$, or
$d\alpha^m\ge k\ge d\alpha^{m+1}$ for some $0\le m\le s-1$. In the first case we
obtain
\[
p=c\alpha^r=c\alpha^{m+r}\le k\alpha^r\le c\alpha^{m+1+r}=p,
\]
which gives $k\alpha^r=p$, so that $k\in C$. A similar argument and conclusion
applies in the alternative case, and so we conclude that
$C=\{k\in X_n: c\le k\le d\}$, as required.
\end{proof}

\begin{rem}
It follows from Proposition~1.5 that the digraph of any $\alpha\in\On$ is a
forest of rooted, labelled trees (each vertex is labelled by a member of $X_n$,
and each component is a tree with a distinguished root vertex in the unique
fixed point of the component). It is easy to verify that this is a necessary and
sufficient condition on a subsemigroup $S$ of $T_n$ to ensure that $S$ is
\emph{aperiodic} or \emph{combinatorial}, meaning that all the subgroups of $S$
are trivial. This is also equivalent to Green's $\mathscr{H}$ relation being
trivial.
\end{rem}

\section{The semigroup of order-preserving mappings}

Let $S_r=\{m_1,m_2,\dots,m_r\}$ with $1\le m_1<m_2<\cdots<m_r\le n$. An
order-preserving mapping $\alpha\colon S_r\to X_n$ is determined by the choice
of $p_1,p_2,\dots,p_r$ such that $1\le p_1\le p_2\le\cdots\le p_r\le n$, where
$m_i\alpha=p_i$ $(i=1,2,\dots,r)$. We can then represent $\alpha$ as a sequence
of $r$ $A$'s and $n$ $B$'s, consisting of $p_1-1$ $B$'s followed by an $A$;
$p_2-p_1$ $B$'s followed by an $A$; \dots; $p_r-p_{r-1}$ $B$'s followed by an
$A$; and finally $n-p_r+1$ $B$'s.

For example, if $n=6$ and $S_r=\{1,2,4,6\}$ then the map
\[
\begin{pmatrix} 1 & 2 & 4 & 6\\ 1 & 3 & 3 & 5\end{pmatrix}
\]
is associated with the sequence $ABBAABBABB$. Any such sequence is determined by
the choice of the $r$ places among the $r+n-1$ available for the positioning of
the $A$'s. Hence we conclude that the number of order-preserving maps from $S_r$
into $X_n$ is $C(r+n-1,r)=C(r+n-1,n-1)$.

Let $Y_n(\alpha)$ stand for the random variable, the value of which is the rank
of $\alpha$ (that is, $|\im\alpha|$ where $\alpha$ is a randomly selected member
of $\On$).

\begin{numthm}{Lemma}{2.1}
$\E(Y_n)=n^2/(2n-1)$, \quad
$\sigma_{Y_n}=\sqrt{(n-1)/2}\cdot n/(2n-1)$.
\end{numthm}

\begin{proof}
We first show that the number of order-preserving mappings
$\alpha\colon S_r\to X_n$ such that $\rank\alpha=k$ is given by
$C(n,k)\cdot C(r-1,k-1)$ $(1\le k\le r)$: the first factor counts the choices
for $\im\alpha=\{i_1,i_2,\dots,i_k\}$ say, while the second counts the choices
for the set $P_\alpha=\{p_1,p_2,\dots,p_k\}$, where $p_j$ is the greatest member
of $S_r$ such that $p_j\alpha=i_j$ $(1\le j\le k)$. (Note that $p_k$ must equal
$m_r$, and thus $p_j<m_r$ for all $j<k$.) Since $\alpha$ is determined by the
choice of the pair $(\im\alpha,P_\alpha)$ and different pairs determine
different mappings, the assertion is established.

Now take $\alpha\in\On$ with $\rank\alpha=r$, while $\beta$ is a randomly chosen
member of $\On$. Denote the random variable with value $\rank\alpha\beta$ by
$Y_{n,r}$. From the foregoing we obtain
\begin{equation}\tag{2.2}
\Pr(Y_{n,r}=k)=\frac{C(n,k)\,C(r-1,k-1)}{C(r+n-1,n-1)}.
\end{equation}
(Taking $r=n$ thus gives the distribution of $Y_n$.) Hence
\begin{align*}
\E(Y_{n,r}) &= \bigl(C(r+n-1,n-1)\bigr)^{-1}\sum_{k=1}^{r} k\,C(n,k)\,C(r-1,k-1)\\
&= \frac{n}{C(r+n-1,n-1)}\sum_{k=0}^{r-1} C(n-1,k)\,C(r-1,k).
\end{align*}
By Result~1.1(i) this equals
\begin{equation}\tag{2.3}
\frac{n}{C(r+n-1,n-1)}\cdot C(r+n-2,n-1)=\frac{nr}{n+r-1}.
\end{equation}
If one now takes $r=n$, (2.3) yields the stated result for
$\E(Y_n)=\E(Y_{n,n})$. Again by using (2.2) one obtains
\begin{align*}
\E\bigl(Y_{n,r}(Y_{n,r}-1)\bigr)
&= C(r+n-1,n-1)^{-1}\sum_{k=2}^{r} k(k-1)\,C(n,k)\,C(r-1,k-1)\\
&= \frac{n(r-1)}{C(r+n-1,n-1)}\sum_{k=1}^{r-1} C(n-1,k)\,C(r-2,k-1).
\end{align*}
We re-write this using the identity
\[
C(r-2,k-1)=C(r-1,k)-C(r-2,k).
\]
We obtain, using Result~1.1(i),
\begin{align*}
\E\bigl(Y_{n,r}(Y_{n,r}-1)\bigr)
&= \frac{n(r-1)}{C(r+n-1,n-1)}\Bigl(\sum_{k=1}^{r-1}C(n-1,k)\,C(r-1,k)\\
&\qquad\qquad\qquad\qquad\quad -\sum_{k=1}^{r-1}C(n-1,k)\,C(r-2,k)\Bigr)\\
&= \frac{n(r-1)\bigl(C(n+r-2,r-1)-1-(C(n+r-3,r-2)-1)\bigr)}{C(r+n-1,n-1)}\\
&= \frac{n(n-1)\,r(r-1)}{(n+r-1)(n+r-2)}.
\end{align*}
Then from the general fact that
$\sigma_X^2=\E(X(X-1))-\E(X)\bigl(\E(X)-1\bigr)$ we obtain
\begin{equation}\tag{2.4}
\Var(Y_{n,r})=\frac{n(n-1)\,r(r-1)}{(n+r-1)^2(n+r-2)}
\end{equation}
and in particular $\Var(Y_{n,r})<n$ for all values of $r=1,2,\dots,n$. The value
of the standard deviation of $Y_n$ is obtained by putting $r=n$ and taking
square roots, thus completing the proof.
\end{proof}

\begin{numthm}{Theorem}{2.5}
Let $Z_k$ be the random variable with value
$n^{-1}\rank(\alpha_1\alpha_2\cdots\alpha_k)$, where
$\alpha_1,\alpha_2,\dots,\alpha_k$ are randomly selected members of $\On$. Then
\[
M_k=\lim_{n\to\infty}\E(Z_k)=1/(1+k)\quad\text{for all } k\ge1.
\]
\end{numthm}

\begin{rem}
If we replace $\On$ by $T_n$ in this result then we obtain $M_k\sim 2/k$: see
\cite{H88}.
\end{rem}

\begin{proof}
We prove by induction on $k$ both the stated value for $M_k$ and that
$\sigma_{Z_k}\to0$ as $n\to\infty$; that this is true for $k=1$ follows from
Lemma~2.1. Assume that our claim holds for some arbitrary value $k-1\ge1$ and let
$\delta,\epsilon>0$ be given. From Chebyshev's Inequality, for any positive
number $R$ we have
\[
\Pr\bigl(|Z_{k-1}-\mu_{k-1}|\ge R\sigma_{Z_{k-1}}\bigr)\le 1/R^2,
\quad\text{where}\quad \mu_{k-1}=\E(Z_{k-1}).
\]
Choose $R$ so that $1/R^2<\epsilon$. Then choose $N$ sufficiently large so that
$R\sigma_{Z_{k-1}}<\delta/2$ and $|\mu_{k-1}-M_{k-1}|<\delta/2$ for all $n\ge N$.
Thus, for all $n\ge N$, $|Z_{k-1}-M_{k-1}|\ge\delta$ implies that
$|Z_{k-1}-\mu_{k-1}|+|\mu_{k-1}-M_{k-1}|\ge\delta$, and hence
$|Z_{k-1}-\mu_{k-1}|\ge\delta/2$. Therefore
\[
\Pr\bigl(|Z_{k-1}-M_{k-1}|\ge\delta\bigr)
\le\Pr\bigl(|Z_{k-1}-\mu_{k-1}|\ge\delta/2\bigr)<\epsilon
\]
for all $n\ge N$. We conclude that
$\Pr(|Z_{k-1}-M_{k-1}|<\delta)\ge1-\epsilon$ for all sufficiently large values
of $n$. Now, by (2.3),
\[
\E(Z_k\mid Z_{k-1}=r/n)=\frac{r}{n+r-1}=\frac{r/n}{1+r/n-1/n},
\]
and this is increasing in $r$. Hence
\begin{align*}
\E(Z_k) &= \sum_{r=1}^{n}\Pr(Z_{k-1}=r/n)\,\E(Z_k\mid Z_{k-1}=r/n)\\
&\ge \sum_{r/n\ge M_{k-1}-\delta}\Pr(Z_{k-1}=r/n)\,
      \E(Z_k\mid Z_{k-1}=M_{k-1}-\delta), \tag{2.6}
\end{align*}
so that
\begin{equation}\tag{2.7}
\E(Z_k)\ge\frac{n(1-\epsilon)(1/k-\delta)}{n+n(1/k-\delta)-1}
\to\frac{(1-\epsilon)(1/k-\delta)}{1+1/k-\delta}
=\frac{(1-\epsilon)(1-k\delta)}{1+k-k\delta}
\quad\text{as } n\to\infty,
\end{equation}
which approaches $1/(1+k)$ as $\delta,\epsilon\to0$. Similarly, splitting the
sum in (2.6) into two sums corresponding to $r/n\le M_{k-1}+\delta$ and
$r/n>M_{k-1}+\delta$, and taking appropriate bounds yields
\begin{equation}\tag{2.8}
\E(Z_k)\le\frac{n(1/k+\delta)}{n+n(1/k+\delta)-1}+\epsilon\cdot1
\to\frac{1+k\delta}{1+k+k\delta}+\epsilon
\end{equation}
which also approaches $1/(1+k)$ as $\delta,\epsilon\to0$ and (2.7) and (2.8)
together yield that $M_k=1/(1+k)$. We complete the proof by verifying
inductively that $\sigma_{Z_k}\to0$ as $n\to\infty$.

Write $W_k=nZ_k$, $p_i=\Pr(W_k=i)$, $p_{ij}=\Pr(W_k=i\mid W_{k-1}=j)$,
$q_j=\Pr(W_{k-1}=j)$, and $\mu_j=\E(W_k\mid W_{k-1}=j)$ for $1\le i,j\le n$.
For the moment let $a$ be arbitrary and consider
\[
\sum_{i=1}^{n}p_i(i-a)^2
=\sum_{i=1}^{n}\Bigl(\sum_{j=1}^{n}p_{ij}q_j\Bigr)(i-a)^2
=\sum_{j=1}^{n}q_j\Bigl(\sum_{i=1}^{n}p_{ij}(i-a)^2\Bigr).
\]
Write $i-a$ as $i-\mu_j+\mu_j-a$ and expand to obtain
\begin{equation}\tag{2.9}
\sum_{j=1}^{n}q_j\Bigl(\sum_{i=1}^{n}p_{ij}(i-\mu_j)^2
+2\sum_{i=1}^{n}p_{ij}(i-\mu_j)(\mu_j-a)
+\sum_{i=1}^{n}p_{ij}(a-\mu_j)^2\Bigr).
\end{equation}
Now $\sum_{i=1}^{n}p_{ij}(i-\mu_j)^2=\sigma_{Y_{n,j}}^2=O(n)$ by (2.4). Also
$\sum_{i=1}^{n}p_{ij}(i-\mu_j)$ is identically zero, so that the sum of the
first two terms in (2.9) is $O(n)$. The remaining term equals
\begin{equation}\tag{2.10}
\sum_{j=1}^{n}\Bigl(q_j(a-\mu_j)^2\sum_{i=1}^{n}p_{ij}\Bigr)
=\sum_{j=1}^{n}q_j(a-\mu_j)^2.
\end{equation}
By Chebyshev's Inequality, for any positive value of $R$,
\[
\Pr\bigl(|W_{k-1}-\mu_{k-1}|\ge R\sigma_{W_{k-1}}\bigr)\le1/R^2.
\]
Let $\epsilon>0$ be arbitrary, and take $R$ such that $1/R^2<\epsilon$. Since
$R\sigma_{W_{k-1}}=o(n)$, it follows that there exist intervals $I_n$ of length
$o(n)$ such that the probability that $W_{k-1}$ lies outside $I_n$ is less than
$\epsilon$. Take $i=i(n)$ to be a fixed member of $I_n\cap X_n$, and put
$a=\E(W_k\mid W_{k-1}=i)$. Now recall that
\[
\E(W_k\mid W_{k-1}=r)=\frac{rn}{r+n-1}=f(r)\ \text{say},\qquad
|f'(r)|=\frac{n(n-1)}{(n+r-1)^2}<1\quad\text{for all } r\ge1.
\]
Hence $|a-\E(W_k\mid W_{k-1}=j)|\le|i-j|=o(n)$ for all $j\in I_n\cap X_n$.
Therefore the summation $\Sigma$ of (2.10) is such that
\begin{align*}
\Sigma &= \sum_{j\in I_n}\Pr(W_{k-1}=j)\bigl(a-\E(W_k\mid W_{k-1}=j)\bigr)^2\\
&\quad +\sum_{j\notin I_n}\Pr(W_{k-1}=j)\bigl(a-\E(W_k\mid W_{k-1}=j)\bigr)^2\\
&\le o(n^2)+\epsilon\cdot n^2.
\end{align*}
Therefore we obtain that
\[
\sigma_{W_k}^2\le\sum_{i=1}^{n}p_i(i-a)^2\le O(n)+o(n^2)+\epsilon\cdot n^2,
\]
and hence that $\sigma_{Z_k}^2\le o(1)+\epsilon$, and since $\epsilon$ was
arbitrary, it follows that $\sigma_{Z_k}^2=o(1)$, whence $\sigma_{Z_k}=o(1)$, as
required.
\end{proof}

Let $Z_n(\alpha)$ now denote the random variable the value of which is the number
of fixed points of a randomly selected $\alpha\in\On$. By Proposition~1.5,
$Z_n(\alpha)$ also equals the component number of the digraph of $\alpha$, and
also the number of cycle points.

The number of mappings $\alpha\in\On$ with $k\alpha=k$ for all $k\in X_n$ is the
product of the number of order-preserving mappings from $X_{k-1}$ to $X_k$ and
the number of such maps from $X_{n-k}$ to $X_{n-k+1}$. By the remark preceding
the proof of Lemma~2.1 the number of order-preserving maps from a set of order
$r$ to a set of order $n$ is $C(r+n-1,n-1)$, whence we conclude that
\[
\Pr(k\alpha=k)=\frac{C(2k-2,k-1)\,C(2n-2k,n-k)}{C(2n-1,n-1)}.
\]
Hence
\begin{align*}
\E(Z_n) &= \sum_{k=1}^{n}\Pr(k\alpha=k)\\
&= C(2n-1,n-1)^{-1}\sum_{k=1}^{n}C(2k-2,k-1)\,C(2n-2k,n-k)\\
&= C(2n-1,n-1)^{-1}\sum_{k=0}^{m}C(2k,k)\,C(2m-2k,m-k),
\end{align*}
where $m=n-1$. From Result~1.1(ii) one now obtains the first statement in the
next result.

\begin{numthm}{Theorem}{2.11}
$\E(Z_n)=4^{n-1}/C(2n-1,n-1)$. Moreover
$\bigl(\E(Z_n)-\sqrt{n\pi}/2\bigr)\to0$ as $n\to\infty$.
\end{numthm}

\begin{proof}
To verify that the mean is of order $\sqrt{n\pi}/2$ one readily shows that
\[
\frac{2(\E(Z_n))^2}{n}
=\frac{(2\cdot4\cdot6\cdots(2n-2))^2\cdot2n}{(1\cdot3\cdot5\cdots(2n-1))^2}.
\]
This is the product of Wallis which approaches $\pi/2$, from which we infer that
$\E(Z_n)\sim\sqrt{n\pi}/2$. In the course of the proof of the Wallis limit it is
established that
\[
\frac{(2\cdot4\cdot6\cdots2n)^2}{(1\cdot3\cdot5\cdots(2n-1))^2(2n+1)}
\le\frac{\pi}{2}\le
\frac{(2\cdot4\cdot6\cdots2n)^2}{(1\cdot3\cdot5\cdots(2n-1))^2\,2n}
\]
from which we obtain that
\[
\frac{4\E^2(Z_n)}{2n+1}\le\frac{\pi}{2}\le\frac{4\E^2(Z_n)}{2n}.
\]
This yields
\[
0\le4\E^2(Z_n)-n\pi\le4\E^2(Z_n)/(2n+1).
\]
We deduce from above that the ratio on the right approaches $\pi/2$, and thus
\[
4\E^2(Z_n)-n\pi=O(1),
\]
from which it follows that $\E(Z_n)-\sqrt{n\pi}/2\to0$.
\end{proof}

\begin{rem}[Remarks]
The closest analogue of this result for $T_n$ concerns the random variable
$X_n(\alpha)$ $(\alpha\in T_n)$, the value of which is the order of the stable
range of $\alpha$ (which corresponds to the number of cycle points of the
digraph of $\alpha$). It is shown in \cite{HW} that
\[
\bigl(\sqrt{n\pi/2}-\E(X_n)\bigr)\to\tfrac13.
\]
The proof of this result is based on Stirling's factorial approximation formula
and a result of Ramanujan concerning the `first half' of the McLaurin series for
$e^n$. By contrast, the component number of $\alpha\in T_n$ has a mean that only
increases logarithmically with $n$. It is shown in \cite{HW} that if
$C_n(\alpha)$ denotes the number of components of a randomly selected
$\alpha\in T_n$ then
\[
\E(C_n)-\tfrac12\log n\to\tfrac12(\gamma+\log2)
\]
where $\gamma$ is the Euler--Mascheroni constant. The first proof of this fact is
due to Kruskal \cite{Kruskal}. From Lemma~2.1 it follows readily that the mean
value of $Y_n(\alpha)$ (the rank of a random $\alpha\in\On$) satisfies
\[
\bigl(\E(Y_n)-n/2\bigr)\downarrow\tfrac14.
\]
The corresponding result for $R_n(\alpha)$, the rank of a random
$\alpha\in T_n$, is noted in \cite{HW}:
\[
\E(R_n)-n(1-e^{-1})\to e^{-1}/2.
\]
\end{rem}

Katz \cite{Katz} considered the problem of finding the probability of
`indecomposability' of a randomly selected $\alpha\in T_n$, by which was meant
the probability that $\alpha$ has only one component. He proved that this
probability is of order $\sqrt{\pi/2n}$. We shall now solve the corresponding
problem for $\On$. The `ballot' argument involved also furnishes us with the
value of $|\Cn|$.

Recall from the first section the natural one-to-one correspondence between
$\alpha\in\On$ and a ballot in which candidate $A$ polls $n$ votes while
candidate $B$ polls $n-1$ votes. Observe that $i\in X_n$ is a fixed point of
$\alpha\in\On$ if and only if in the count of the corresponding ballot the $i$th
$A$ appears between the $(i-1)$th and the $i$th $B$. In terms of the count, this
corresponds to a tie after $2i-2$ votes, with $A$ polling the next vote. Now to
say that $\alpha\in\On$ is indecomposable means exactly that $\alpha$ has a
unique fixed point, which, formulated in terms of the corresponding ballot
count, means precisely that $A$ never surrenders the lead once he has taken it
(i.e.\ once $A$ is strictly ahead in the count, he stays strictly ahead).

We consider a more general ballot in which $A$ polls $n+r$ votes while his
opponent $B$ polls $n$ votes $(r\ge1)$, and we seek the probability that $A$
never relinquishes the lead once he has gained it. Let $Y_1$ denote the set of
ballots of this type, let $Y_2$ denote the set of ballots in which $A$ never
trails at any point in the count, and let $Y=Y_1\cup Y_2$. Consider
$\alpha\in Y$, and let $(p,0)$ be the \emph{final} point on the lattice path of
$\alpha$ where the count is tied. Let $\alpha'$ be the count of the lattice path
which results from reflection of that portion of the path of $\alpha$ from
$(0,0)$ to $(p,0)$ in the $x$-axis. Then $\alpha'\in Y$ and we observe that the
mapping ${}'\colon Y\to Y$ defines a bijection between $Y_1$ and $Y_2$ with
$(\alpha')'=\alpha$. In fact, ${}'$ defines a bijection between
$Y_1\setminus Y_2$ and $Y_2\setminus Y_1$ as $Y_1\cap Y_2$ consists precisely of
all those counts where $A$ leads throughout, which coincides with the set of
fixed points of the mapping ${}'$.

We next find $|Y_2|\ (=|Y_1|)$ by calculating the order of the complementary set
of counts, which consists of all those counts during which $A$ trails $B$ at some
point. The number of paths that correspond to counts where $A$ falls behind at
some point equals the number of lattice paths from $(0,0)$ to $(2n+r,r)$ which
touch or cross the line $y=-1$. By the Reflection Principle (Result~1.3) this
equals the number of lattice paths from $(-2,0)$ to $(2n+r,r)$, which is
$C(2n+r,n-1)$; the $n-1$ corresponds to the $n-1$ places where the path dips,
leaving $n+r+1$ rises. Hence the probability that $A$ trails $B$ at some time
during the count is $C(2n+r,n-1)/C(2n+r,n)=n/(n+r+1)$. Hence the probability of
the complementary event that $A$ never trails is
$1-n/(n+r+1)=(r+1)/(n+r+1)$, and this equals the probability that $A$ never
loses his lead once he has gained it. By taking $r=1$ and replacing $n$ by $n-1$
we obtain our required result.

\begin{numthm}{Theorem}{2.12}
The probability that $\alpha\in\On$ is indecomposable is $2/(n+1)$.
\end{numthm}

\section{The semigroup of all decreasing and order-preserving mappings}

We saw in the preceding section that a mapping $\alpha\in\On$ is indecomposable
if and only if in the corresponding ballot count the winning candidate $A$ never
relinquishes his lead once he has gained it, and that this collection of counts
is in bijection with the set of counts where $A$ never trails his opponent. We
now ask what type of mapping this latter set of counts represents. Candidate $A$
never trails in the count if and only if the $i$th vote for $A$ appears before
the $i$th vote for $B$ for $i=1,2,\dots,n$; in other words, the number of votes
for $A$ preceding the $i$th for $B$ in the sequence is always greater than or
equal to $i$. In the corresponding mapping, since each $i\in X_n$ is mapped to
the number of the next $B$ following the $i$th $A$ in the sequence, this
condition is equivalent to $i\alpha\le i$ for all $i\in X_n$. Therefore the
collection of mappings corresponding to the set of counts in which the winning
candidate never trails is $\Cn$, the semigroup of order-preserving and
decreasing mappings, and $|\Cn|$ equals the cardinality of the set of all
indecomposable mappings in $\On$. This yields our next result.

\begin{numthm}{Theorem}{3.1}
$|\Cn|=C_n$, the $n$th Catalan number.
\end{numthm}

\begin{proof}
From the preceding argument, together with Theorem~2.12 and Result~1.2, we have
\[
|\Cn|=|\{\alpha\in\On:\alpha\text{ indecomposable}\}|
=\frac{2}{n+1}\binom{2n-1}{n-1}=\frac{1}{n+1}\binom{2n}{n},
\]
which is the $n$th Catalan number.
\end{proof}

\begin{rem}[Remarks]
For each collection of counts of a given ballot there corresponds a dual
collection obtained by reversing the count sequence of the ballot. Given a count
where $A$ never relinquishes his lead once he has gained it, the dual count has
the property that $A$ never relinquishes his eventual winning lead once he has
taken it. In our case, since the winning margin for $A$ is one, the
characteristic conditions of both types of count are identical, so that the set
of indecomposable mappings of $\On$ is self-dual under this correspondence. On
the other hand, the dual of the condition that $A$ never trails in the count is
that $A$'s lead never exceeds his eventually winning margin, which in our case is
tantamount to saying that $A$'s lead never exceeds one. This is to say that the
number of the next $B$ to follow the $i$th vote for $A$ in the count is at least
$i$, or in terms of the corresponding mapping $\alpha$, $i\alpha\ge i$ for all
$i\in X_n$. Hence the dual to $\Cn$ under this correspondence is the isomorphic
semigroup of all increasing and order-preserving mappings on $X_n$.
\end{rem}

There are natural bijections between $\Cn$ and other standard collections of
order $C_n$ which thus yield alternative proofs of Theorem~3.1. We give the
specific connection between $\Cn$ and the set of all binary trees on $n$ source
nodes. In this context it is more natural to take the base set $X_n$ to be
$\{0,1,2,\dots,n-1\}$.

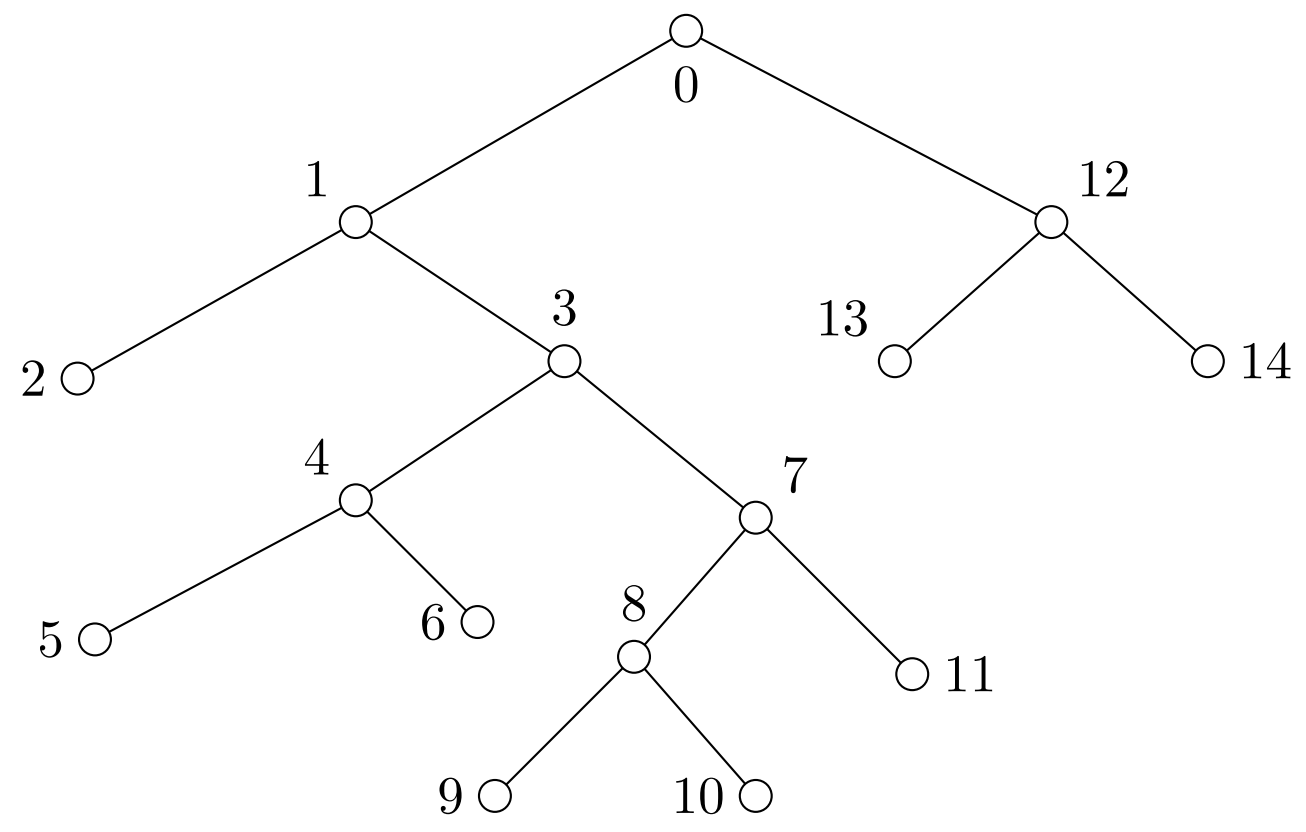
\begin{figure}[ht]
\centering
\begin{tikzpicture}[scale=1.15,
    nd/.style={circle,draw=black,fill=white,inner sep=0pt,minimum size=6pt}]
  \coordinate (n0)  at (4.3,5.0);
  \coordinate (n1)  at (2.4,3.9);
  \coordinate (n12) at (6.4,3.9);
  \coordinate (n2)  at (0.8,3.0);
  \coordinate (n3)  at (3.6,3.1);
  \coordinate (n13) at (5.5,3.1);
  \coordinate (n14) at (7.3,3.1);
  \coordinate (n4)  at (2.4,2.3);
  \coordinate (n7)  at (4.7,2.2);
  \coordinate (n5)  at (0.9,1.5);
  \coordinate (n6)  at (3.1,1.6);
  \coordinate (n8)  at (4.0,1.4);
  \coordinate (n11) at (5.6,1.3);
  \coordinate (n9)  at (3.2,0.6);
  \coordinate (n10) at (4.7,0.6);
  \draw (n0)--(n1) (n0)--(n12)
        (n1)--(n2) (n1)--(n3)
        (n3)--(n4) (n3)--(n7)
        (n4)--(n5) (n4)--(n6)
        (n7)--(n8) (n7)--(n11)
        (n8)--(n9) (n8)--(n10)
        (n12)--(n13) (n12)--(n14);
  \foreach \v in {0,...,14} \node[nd] at (n\v) {};
  \node[below=3pt]       at (n0)  {\small 0};
  \node[above left=1pt]  at (n1)  {\small 1};
  \node[left=2pt]        at (n2)  {\small 2};
  \node[above=3pt]       at (n3)  {\small 3};
  \node[above left=1pt]  at (n4)  {\small 4};
  \node[left=2pt]        at (n5)  {\small 5};
  \node[left=2pt]        at (n6)  {\small 6};
  \node[above right=1pt] at (n7)  {\small 7};
  \node[above=3pt]       at (n8)  {\small 8};
  \node[left=2pt]        at (n9)  {\small 9};
  \node[left=2pt]        at (n10) {\small 10};
  \node[right=2pt]       at (n11) {\small 11};
  \node[above right=1pt] at (n12) {\small 12};
  \node[above left=1pt]  at (n13) {\small 13};
  \node[right=2pt]       at (n14) {\small 14};
\end{tikzpicture}
\caption{A binary tree with seven source nodes.}
\label{fig:tree}
\end{figure}

By a \emph{binary tree} we shall mean a complete binary tree where every node $u$
is either a \emph{source node}, meaning that there are exactly two arcs from $u$
leading to its two sons, or that $u$ is an endpoint or \emph{leaf}. Every node
except the uppermost root of the tree is the son of exactly one source node, so
that it follows that a binary tree with $n$ source nodes has $2n+1$ nodes in all.
We shall also assume that our binary trees are ordered, in the sense that the
nodes are consecutively numbered in the order in which they are met while
carrying out an anti-clockwise search of the tree, beginning at the root. An
example of a binary tree with seven source nodes is given in
Figure~\ref{fig:tree}.

Given a binary tree $B$ on $n$ source nodes (and so with $n+1$ leaf nodes) we
associate a mapping $\alpha\in\Cn$ by the rule that $i\mapsto t$, where $t$ leaf
nodes precede the $i$th source node in the ordering of the nodes of $B$. To see
that $\alpha\in\Cn$ observe that our rule certainly defines an order-preserving
map on $X_n$, and since the final source node precedes its two sons in the
ordering of the nodes, the maximum possible value of $(n-1)\alpha$ is
$n+1-2=n-1$. We can then prove by induction on $n$ that $\alpha$ also has the
decreasing property: remove the two sons of the final source node from $B$ to
give a binary tree $B'$ on $n-1$ source nodes, which, by induction, we may assume
is associated with a member of $\mathcal{C}_{n-1}$ by the above rule. Indeed,
the mapping $\alpha'$ associated with $B'$ is just the restriction of the mapping
$\alpha$ to $X_{n-1}$, and so, by the inductive hypothesis, we have that
$\alpha|X_{n-1}$ is decreasing, and we have already observed that
$(n-1)\alpha\le n-1$, allowing us to conclude that $\alpha\in\Cn$. The mapping
$\alpha\in\mathcal{C}_7$ associated with our binary tree above has its ordered
list of images of the members of $X_7$ the tuple $(0,0,1,1,3,3,6)$.

Conversely, we can associate with any $\alpha\in\Cn$ a binary tree $B$ on $n$
source nodes using the inverse of the foregoing procedure. Suppose inductively
that we have associated with $\alpha'=\alpha|X_{n-1}$ a binary tree $B'$ on $n-1$
source nodes in such a way that $B'$ is associated with $\alpha'$ by the
procedure of the preceding paragraph. Let $k$ denote $(n-1)\alpha$. The tree $B'$
has $n$ leaf nodes. Take the $(k+1)$th leaf node (with $k+1\le n$), and replace
it by a source node to form the binary tree $B$. By construction, the mapping
associated with $B$ is $\alpha$, as required.

The number of indecomposable $\alpha\in\Cn$ is easily obtained through our ballot
correspondence. In terms of the corresponding count,
\begin{enumerate}
\item[(i)] $\alpha\in\Cn$ if and only if $A$ never trails $B$ in the count,
\item[(ii)] $\alpha$ is indecomposable if and only if $A$ never relinquishes the
lead once gained.
\end{enumerate}
The conjunction of this pair of conditions is equivalent to the single
constraint that $A$ leads the count throughout. The enumeration of such counts is
the classical Bertrand--Whitworth ballot problem, readily solved by the
Reflection Principle (see \cite{Feller} or \cite{HP}), and in general, the
probability of such a ballot is $r/(2n+r)$. Again, upon putting $r=1$ and
replacing $n$ by $n-1$ we see that the proportion of such mappings in $\Cn$ is
\[
\frac{1}{2n-1}\binom{2n-1}{n-1}=C_{n-1}.
\]
This yields

\begin{numthm}{Proposition}{3.2}
The proportion of indecomposable mappings in $\Cn$ is
\[
C_{n-1}/C_n=(n+1)/\bigl(2(2n-1)\bigr)\downarrow\tfrac14\quad\text{as } n\to\infty.
\]
\end{numthm}

It is interesting here that the probability decreases to $\frac14$, and not to
$0$, as is the case with $\On$ and $T_n$. This suggests the existence of a
limiting distribution for $Z_n(\alpha)$, where $Z_n(\alpha)$ denotes the
component number of a randomly chosen $\alpha\in\Cn$. Our determination of this
limiting distribution involves the use of the standard Catalan identity for
which we can now offer a new proof.

\begin{proof}[Proof of Result 1.4]
Partition $\Cn$ into $n$ disjoint sets $K_1,K_2,\dots,K_n$, where
$\alpha\in K_k$ if $k$ is the greatest integer such that $k\alpha=k$. Note that
$\alpha\in K_1$ if and only if $1$ is the unique fixed point of $\alpha$, which
is to say that $\alpha$ is indecomposable. Then
\[
|\Cn|=\sum_{k=1}^{n}|K_k|.
\]
Now
$|K_k|=|\mathcal{C}_{k-1}|\,|\{\alpha\in\mathcal{C}_{n-k+1}:\alpha\text{ indecomposable}\}|
=C_{k-1}C_{n-k}$, since, in general, the number of indecomposable mappings in
$\Cn$ is $C_{n-1}$. Hence
\[
|\Cn|=\sum_{k=1}^{n}C_{k-1}C_{n-k}. \qedhere
\]
\end{proof}

We can use this result to show that the expected number of components of a random
$\alpha\in\Cn$ is a monotonic function of $n$ increasing to a limit of $3$. This
contrasts with Theorem~2.11 which shows that the mean component number in $\On$
is of order $\sqrt n$.

\begin{numthm}{Proposition}{3.3}
$\E(Z_n)=3n/(n+2)$ where $Z_n(\alpha)$ is the component number of a random
$\alpha\in\Cn$.
\end{numthm}

\begin{proof}
We have
\begin{align*}
\E(Z_n) &= \sum_{k=1}^{n}\Pr(k\alpha=k)=C_n^{-1}\sum_{k=1}^{n}C_{k-1}C_{n+1-k}\\
&= C_n^{-1}\Bigl(\sum_{k=1}^{n+1}C_{k-1}C_{n+1-k}-C_nC_0\Bigr)
 =(C_{n+1}-C_n)/C_n\\
&= \frac{2(2n+1)}{n+2}-1=\frac{3n}{n+2}\uparrow3. \qedhere
\end{align*}
\end{proof}

Let $N(n,k)=|\{\alpha\in\Cn:\alpha\text{ has exactly } k \text{ fixed points}\}|$.

\begin{numthm}{Lemma}{3.4}
\[
N(n,k)=\frac{k}{2n-k}\binom{2n-k}{n},\quad\text{for } k=1,2,\dots,n.
\]
\end{numthm}

\begin{proof}
Using Proposition~3.2 we have
\begin{align*}
N(n,1) &= |\{\alpha\in\Cn:\alpha\text{ indecomposable}\}|=C_{n-1}\\
&= (1/n)\,C(2n-2,n-1)=\bigl(1/(2n-1)\bigr)\,C(2n-1,n),
\end{align*}
as asserted. Using the argument of the proof of Result~1.4, we focus attention on
the final fixed point of a mapping in $\Cn$ and conclude that
\begin{align*}
N(n,2) &= \sum_{k=1}^{n-1}N(k,1)\,N(n-k,1)=\sum_{k=1}^{n-1}C_{k-1}C_{n-k-1}\\
&= \sum_{k=1}^{n-1}C_{k-1}C_{(n-1)-k}=C_{n-1}
\end{align*}
by Result~1.4, and $C_{n-1}=\bigl(2/(2n-2)\bigr)C(2n-2,n)$, as required.

We next establish the recurrence
\begin{equation}\tag{3.5}
N(n,k+1)=N(n,k)-N(n-1,k-1)\quad\text{for } k=1,2,\dots,n-1.
\end{equation}
The above calculation shows that (3.5) is valid for $k=1$ (as $N(n,0)=0$). Take
$k\ge2$. By considering the final fixed point $t$ of $\alpha\in\Cn$ with exactly
$k+1$ fixed points we see that
\begin{equation}\tag{3.6}
N(n,k+1)=\sum_{t=k+1}^{n}N(t-1,k)\,N(n-t+1,1).
\end{equation}
By considering the penultimate fixed point $t$ we also have
\begin{equation}\tag{3.7}
N(n,k+1)=\sum_{t=k}^{n-1}N(t-1,k-1)\,N(n-t+1,2).
\end{equation}
But we know that $N(n-t+1,2)=N(n-t+1,1)$, and so, by substituting accordingly in
(3.7), we obtain
\[
N(n,k+1)=\sum_{t=k}^{n}N(t-1,k-1)\,N(n-t+1,1)-N(n-1,k-1)\,N(1,1).
\]
Therefore $N(n,k+1)=N(n,k)-N(n-1,k-1)$, as asserted.

It is now a routine matter, using a double induction on $n$ and $k$, to verify
that the solution to our recurrence is given by
\[
N(n,k)=\bigl(k/(2n-k)\bigr)\,C(2n-k,n)\quad\text{for } k=1,2,\dots,n. \qedhere
\]
\end{proof}

\begin{numthm}{Theorem}{3.8}
The probability distribution for the random variable $Z_n(\alpha)$, the component
number of a random $\alpha\in\Cn$, is given by
\begin{equation}\tag{3.9}
\Pr(Z_n=k)=p(n,k)=\frac{k(n+1)(n-1)(n-2)\cdots(n-k+1)}{2(2n-1)(2n-2)\cdots(2n-k)}
\quad\text{for } k=1,2,\dots,n.
\end{equation}
As $n$ increases, $Z_n$ approaches in distribution and in moments the random
variable $Z$, where $p_k=\Pr(Z=k)=k/2^{k+1}$ for $k=1,2,\dots$. In particular,
the limiting values of the mean and variance are $3$ and $2$ respectively.
Furthermore, for $n>1$, the probability sequence is bi-modal, with modes $k=1$
and $k=2$, and is strictly monotonically decreasing for $k>2$, features shared
with the limiting distribution $Z$.
\end{numthm}

\begin{proof}
Except for the statement about moments, this is now immediate: from
Theorem~3.1 and Lemma~3.4 we obtain
\[
p(n,k)=\frac{N(n,k)}{C_n}=\frac{k}{2n-k}\,C(2n-k,n)\Big/
\Bigl(\frac{1}{n+1}\,C(2n,n)\Bigr),
\]
which simplifies to yield the expression (3.9). That
$\lim_{n\to\infty}p(n,k)=k/2^{k+1}$ follows on dividing each of the last $k$
terms in the numerator and denominator by $n$, so that $Z_n$ approaches $Z$ in
distribution. From (3.9) we deduce the recurrence
\begin{equation}\tag{3.10}
p(n,k+1)=p(n,k)\cdot\frac{(k+1)(n-k)}{k(2n-k-1)}\quad\text{for } k=1,2,\dots,n-1
\end{equation}
from which the statements on modes and on monotonicity follow. The statement on
convergence in moments can be deduced from the following fact.

\begin{numthm}{Lemma}{3.11}
If $n>2$ then $p(n,k)>p_k$ for $k=1,2,3$ and $p(n,k)<p_k$ for $k\ge4$.
\end{numthm}

\noindent\textit{Proof.} Direct calculation using (3.9) allows one to verify the
assertions of the Lemma for $k=1,2,3$ and $4$. It remains only to observe that if
$p(n,k)<p_k$ for some $k$, then (3.10) yields
\[
p(n,k+1)<p_{k+1}\,\frac{2(n-k)}{2n-k-1},
\]
from which we infer that $p(n,k+1)<p_{k+1}$. The result now follows by induction
on $k$. \hfill$\square$

\medskip
To show convergence in moments, first observe that $\E(Z^t)$ certainly exists for
$t=1,2,\dots$ (indeed the values of $\E(Z^t)$ are integers, and can be calculated
recursively). Let $t$ be a positive integer, and let $\epsilon>0$ be given.
Choose $N\ge3$ such that
\[
\sum_{k=N+1}^{\infty}k^t p_k<\epsilon/3.
\]
Then, for any positive integer $n$,
\begin{align*}
|\E(Z^t)-\E(Z_n^t)|
&\le\Bigl|\sum_{k=1}^{N}k^t\bigl(p_k-p(n,k)\bigr)\Bigr|
   +\sum_{k=N+1}^{\infty}k^t\bigl(p_k-p(n,k)\bigr)\\
&\le\Bigl|\sum_{k=1}^{N}k^t\bigl(p_k-p(n,k)\bigr)\Bigr|
   +\sum_{k=N+1}^{\infty}k^t p_k+\sum_{k=N+1}^{\infty}k^t p(n,k).
\end{align*}
Since $Z_n$ approaches $Z$ in distribution, it follows that for sufficiently
large $n$ the first term is less than $\epsilon/3$; the same is true of the
second term, and of the third which is bounded by the second. Hence, for all $n$
sufficiently large, $|\E(Z^t)-\E(Z_n^t)|<\epsilon$, and so $Z_n$ approaches $Z$
in moments also.
\end{proof}

\begin{rem}
Using the fact that $p(n,k)<p_k$ for all $k$ sufficiently large, it can be shown
that $\E(Z_n^t)<\E(Z^t)$ for all positive integers $n$ and $t$.
\end{rem}

The probability distribution for the component number of a random
$\alpha\in\On$ can be calculated explicitly in a similar fashion. Let us now
write $Z_n(\alpha)$ for the component number of a randomly chosen
$\alpha\in\On$, and let $N(n,k)$ be the cardinality of
$\{\alpha\in\On:\alpha\text{ has exactly } k\text{ components}\}$.

\begin{numthm}{Theorem}{3.12}
Let $p(n,k)=\Pr(Z_n=k)$, where $Z_n$ is the component number of a random
$\alpha\in\On$. Then
\begin{equation}\tag{3.13}
p(n,k)=\frac{2k(n-1)(n-2)\cdots(n-k+1)}{(n+1)(n+2)\cdots(n+k)}
\quad\text{for } k=1,2,\dots,n.
\end{equation}
The distribution of $Z_n$ has a unique mode at
$\bigl[\frac12\bigl(\sqrt{1+2n}-1\bigr)\bigr]$, the greatest integer not exceeding
$\frac12\bigl(\sqrt{1+2n}-1\bigr)$, unless $1+2n$ is a perfect square, in which
case $\bigl[\frac12\bigl(\sqrt{1+2n}+1\bigr)\bigr]$ is another mode. The
distribution decreases monotonically on either side of the mode(s) with the
minimum occurring at $k=n$. Furthermore, $p(n,k)\sim2k/n$ as $n\to\infty$.
\end{numthm}

\begin{proof}
Assertion (3.13) for $k=1$ is Theorem~2.12; in this case $N(n,1)=C_n$. In order
to establish (3.13) for $k>1$ it is enough, by Result~1.2, to verify that
\begin{equation}\tag{3.14}
N(n,k)=\frac{k\,(2n)!\,(n-1)(n-2)\cdots(n-k+1)}{n!\,(n+k)!}
\quad\text{for } k=1,2,\dots,n.
\end{equation}
We do this by proving the analogue of Lemma~3.4.

\begin{numthm}{Lemma}{3.15}
$N(n,k+1)=N(n+1,k)-2N(n,k)-N(n,k-1)$ for $k=1,2,\dots,n-1$.
\end{numthm}

\noindent\textit{Proof.} Let $\alpha$ be a typical member of $\On$ with exactly
two components, and thus two fixed points, $p<q$, say. Bearing in mind the
convexity of the orbits of $\alpha$ (Proposition~1.5), and letting $t$ be the
maximum of the orbit of $p$, we see that
\[
N(n,2)=\sum_{t=1}^{n-1}N(t,1)\,N(n-t,1)=\sum_{t=1}^{n-1}C_tC_{n-t}
=\sum_{t=2}^{n}C_{t-1}C_{n+1-t}.
\]
Applying Result~1.4 then yields that
\[
N(n,2)=C_{n+1}-C_0C_n-C_nC_0=C_{n+1}-2C_n.
\]
In other words, $N(n,2)=N(n+1,1)-2N(n,1)$, in agreement with the statement of
our lemma.

Next, for $k\ge2$, by considering the final, and then the penultimate fixed point
of a mapping from $\On$ with exactly $k+1$ fixed points, we obtain the
expressions
\begin{align}
N(n,k+1) &= \sum_{t=k}^{n-1}N(t,k)\,N(n-t,1), \tag{3.16}\\
N(n,k+1) &= \sum_{t=k-1}^{n-2}N(t,k-1)\,N(n-t,2). \tag{3.17}
\end{align}
Replacing $N(n-t,2)$ by $N(n-t+1,1)-2N(n-t,1)$ in (3.17) yields that
\[
N(n,k+1)=\sum_{t=k-1}^{n-2}N(t,k-1)\,N(n-t+1,1)
-2\sum_{t=k-1}^{n-2}N(t,k-1)\,N(n-t,1).
\]
From (3.16) we see that the first term can be written as
\[
N(n+1,k)-N(n-1,k-1)\,N(2,1)-N(n,k-1)\,N(1,1);
\]
while the summation in the second term equals
\[
N(n,k)-N(n-1,k-1)\,N(1,1).
\]
Combining these expressions yields
\[
N(n,k+1)=N(n+1,k)-2N(n,k)-N(n,k-1)
\]
as required. \hfill$\square$

\medskip
A routine calculation now shows that the expression for $N(n,k)$ given in (3.14)
is correct, and thus (3.13) is also.

From (3.13) comes the recurrence relation
\begin{equation}\tag{3.18}
p(n,k+1)=p(n,k)\,\frac{k+1}{k}\,\frac{(n-k)}{(n+k+1)}.
\end{equation}
From (3.18) it follows that
\begin{equation}\tag{3.19}
p(n,k+1)>p(n,k)\quad\text{if and only if}\quad
k<\tfrac12\bigl(\sqrt{1+2n}-1\bigr)
\end{equation}
with equality occurring in both sides simultaneously. Using $[x]$ to denote the
greatest integer not exceeding a real number $x$, we see that if
$\frac12(\sqrt{1+2n}-1)$ is not an integer then $Z_n$ has a unique mode at
$[\frac12(\sqrt{1+2n}-1)]$, while otherwise $Z_n$ is bi-modal with modes
$[\frac12(\sqrt{1+2n}\pm1)]$. That the distribution of $Z_n$ falls away
monotonically on either side of the mode(s) follows from (3.19), and that the
minimum value of $p(n,k)$ occurs at $k=n$ follows by comparing
$p(n,1)=2/(n+1)$ with $p(n,n)=1/C(2n-1,n-1)$. Finally, from (3.13),
\[
p(n,k)=\frac{2k}{n}\cdot\frac{n(n-1)\cdots(n-k+1)}{(n+1)(n+2)\cdots(n+k)}
\]
and since the bracketed term approaches $1$ for large $n$, it follows that
$p(n,k)\sim2k/n$ as $n\to\infty$.
\end{proof}

In \cite{Howie}, the cardinality of $E(\On)$, the set of idempotents of $\On$,
was shown to be $F_{2n}$, the $2n$th Fibonacci number. This was achieved by
solving a certain recurrence relation. Determining the order of $E(\Cn)$ is,
however, a simple matter, for any $\alpha\in E(\Cn)$ is determined by
$\im\alpha$, the image of $\alpha$. To see this, first recall that
$\alpha\in T_n$ is an idempotent if and only if $\im\alpha$ and the set of fixed
points of $\alpha$ coincide. If we further suppose that $\alpha\in E(\Cn)$ then,
for any $i\in X_n$, $i\alpha\le i$, $i\alpha$ is a fixed point of $\alpha$, and
since $\alpha$ preserves order, $i\alpha$ must be the greatest fixed point less
than or equal to $i$. Hence an idempotent of $\Cn$ is determined by its image and
thus $|E(\Cn)|=|\{\im\alpha:\alpha\in\Cn\}|$. Now the subsets of $X_n$ which
occur as images of mappings $\alpha\in\Cn$ are exactly the subsets of $X_n$ that
contain $1$. This yields the first statement in our final result.

\begin{numthm}{Theorem}{3.19}
$|E(\Cn)|=2^{n-1}$. Furthermore,
$R=|E(\Cn)|/|\Cn|\sim\sqrt{\pi n}\,(n+1)/2^{n+1}$.
\end{numthm}

\begin{proof}
We have $R=2^{n-1}/C_n$, from which one obtains
$4^{n+1}R^2/\bigl(2n(n+1)^2\bigr)=W_n$, where $W_n$ is the Wallis product as in
the proof of Theorem~2.11. Hence we deduce that
$2^{n+1}R/\bigl(\sqrt{n\pi}\,(n+1)\bigr)\to1$ as $n\to\infty$, as required.
\end{proof}

\begin{rem}
For an investigation of the number of idempotents of finite full transformation
semigroups see \cite{Harris} and \cite{HS}.
\end{rem}



\begin{thebibliography}{99}

\bibitem{Feller} W.~Feller.
\newblock \emph{An Introduction to Probability Theory and its Applications}, 3rd
edn.\ vol.~1 (Wiley, 1968).

\bibitem{Garba} G.~U.~Garba.
\newblock Nilpotents in semigroups of partial one-to-one order-preserving
mappings.
\newblock \emph{Semigroup Forum} \textbf{41} (1991), 1--15.

\bibitem{GH} G.~M.~S.~Gomes and J.~M.~Howie.
\newblock On the ranks of certain semigroups of order-preserving
transformations.
\newblock \emph{Semigroup Forum} \textbf{40} (1990), 1--15.

\bibitem{Harris} B.~Harris.
\newblock A note on the number of idempotents in symmetric semigroups.
\newblock \emph{J.\ Combin.\ Theory Ser.\ A} \textbf{3} (1967), 1234--1235.

\bibitem{HS} B.~Harris and L.~Schoenfeld.
\newblock The number of idempotent elements in symmetric semigroups.
\newblock \emph{J.\ Combin.\ Theory Ser.\ A} \textbf{3} (1967), 122--135.

\bibitem{H88} P.~M.~Higgins.
\newblock The range order of a product of $i$ transformations from a finite full
transformation semigroup.
\newblock \emph{Semigroup Forum} \textbf{37} (1988), 31--36.

\bibitem{H90} P.~M.~Higgins.
\newblock Random products in semigroups of mappings.
\newblock In \emph{Lattices, Semigroups and Universal Algebra} (editors
J.~Almeida et al.), (Plenum Press, 1990), pp.~89--100.

\bibitem{H92} P.~M.~Higgins.
\newblock \emph{Techniques of Semigroup Theory} (Oxford University Press, 1992).

\bibitem{HW} P.~M.~Higgins and E.~J.~Williams.
\newblock Random functions on a finite set.
\newblock \emph{Ars Combin.} \textbf{26A} (1988), 93--102.

\bibitem{HP} P.~Hilton and J.~Pedersen.
\newblock Catalan numbers, their generalization, and their uses.
\newblock \emph{Math.\ Intelligencer} \textbf{13} (1991), 64--75.

\bibitem{Howie} J.~M.~Howie.
\newblock Products of idempotents in certain semigroups of transformations.
\newblock \emph{Proc.\ Roy.\ Soc.\ Edinburgh Sect.\ A} \textbf{17} (1971),
223--236.

\bibitem{Katz} L.~Katz.
\newblock Probability of indecomposability of a random mapping function.
\newblock \emph{Ann.\ Math.\ Statist.} \textbf{26} (1955), 512--517.

\bibitem{Kruskal} M.~D.~Kruskal.
\newblock The expected number of components under a random mapping function.
\newblock \emph{Amer.\ Math.\ Monthly} \textbf{61} (1954), 392--397.

\bibitem{Pin} J.~E.~Pin.
\newblock \emph{Varieties of Formal Languages} (Plenum Press, 1986).

\end{thebibliography}
\end{document}